
\magnification=1200
\input amstex.tex
\documentstyle{amsppt}
\hsize=12.5cm
\vsize=18cm
\hoffset=1cm
\voffset=2cm

\def\DJ{\leavevmode\setbox0=\hbox{D}\kern0pt\rlap
 {\kern.04em\raise.188\ht0\hbox{-}}D}
\footline={\hss{\vbox to 2cm{\vfil\hbox{\rm\folio}}}\hss}
\nopagenumbers
\font\ff=cmr8
\def\txt#1{{\textstyle{#1}}}
\baselineskip=13pt
\def\hf{{\textstyle{1\over2}}}
\def\a{\alpha}
\def\d{{\,\roman d}}
\def\e{\varepsilon}
\def\f{\varphi}
\def\G{\Gamma}
\def\k{\kappa}
\def\s{\sigma}

\def\z{\zeta}
\def\={\;=\;}
\def\zx{\zeta(\hf+ix)}
\def\zt{\zeta(\hf+it)}

\def\no{\noindent}  
\def\R{\Re{\roman e}\,}  \def\s{\sigma}
\def\Z{{\Cal Z}}

\def\no{\noindent}

\def\H{H_j^3({\txt{1\over2}})}
\font\teneufm=eufm10
\font\seveneufm=eufm7
\font\fiveeufm=eufm5
\newfam\eufmfam
\textfont\eufmfam=\teneufm
\scriptfont\eufmfam=\seveneufm
\scriptscriptfont\eufmfam=\fiveeufm
\def\mathfrak#1{{\fam\eufmfam\relax#1}}

\font\tenmsb=msbm10
\font\sevenmsb=msbm7
\font\fivemsb=msbm5
\newfam\msbfam
\textfont\msbfam=\tenmsb
\scriptfont\msbfam=\sevenmsb
\scriptscriptfont\msbfam=\fivemsb
\def\Bbb#1{{\fam\msbfam #1}}

\def \NN {\Bbb N}
\def \CC {\Bbb C}

\def \RR {\Bbb R}
\def \ZZ {\Bbb Z}

\def\rightheadline{{\hfil{\ff
On sums of integrals of powers of the zeta-function in short intervals}
\hfil\tenrm\folio}}

\def\leftheadline{{\tenrm\folio\hfil{\ff
Aleksandar Ivi\'c }\hfil}}
\def\emptyheadline{\hfil}
\headline{\ifnum\pageno=1 \emptyheadline\else
\ifodd\pageno \rightheadline \else \leftheadline\fi\fi}

\topmatter
\title
ON SUMS OF INTEGRALS OF POWERS OF THE ZETA-FUNCTION
IN SHORT INTERVALS
\endtitle
\author
 Aleksandar Ivi\'c
\endauthor
\address{
Aleksandar Ivi\'c, Katedra Matematike RGF-a
Universiteta u Beogradu, \DJ u\v sina 7, 11000 Beograd,
Serbia and Montenegro.
}
\endaddress
\email {\tt aivic\@ivic.matf.bg.ac.yu, ivic\@rgf.bg.ac.yu} \endemail
\abstract
{The modified Mellin transform  $\Z_k(s) = \int_1^\infty |\zx|^{2k}x^{-s}\d x
\;(k\in\NN)$ is used to obtain estimates for
$$
\sum_{r=1}^R\int_{t_r-G}^{t_r+G}|\zt|^{2k}\d t\quad(T < t_1 < \cdots < t_R < 2T),
$$
where $t_{r+1} - t_r \ge G\;(r =1,\ldots, R-1), \,T^\e \le G \le T^{1-\e}$. These
results can be used to derive bounds for the moments of $|\zt|$.}
\endabstract
\endtopmatter

\heading
{\bf 1. Introduction}
\endheading

The (modified) Mellin transforms
$$
\Z_k(s) \;:=\; \int_1^\infty |\zx|^{2k}x^{-s}\d x
\qquad(k \in \NN,\,\s = \R s \ge c(k) > 1),
\leqno(1.1)
$$
where $c(k)$ is such a constant for which the integral in (1.1) converges
absolutely, play an important r\^ole in the theory of
the Riemann zeta-function $\z(s)$ (see [1], [7], [9], [14] and [19] for some of the
relevant works, which contain further references).
The term ``modified" Mellin transform seems appropriate, since
customarily the Mellin transform of $f(x)$ is defined as
$$
 F(s) := \int_0^\infty f(x)x^{s-1}\d x\qquad(s = \s + it\in\CC).
\leqno(1.2)
$$
Note that the lower bound of integration in (1.1) is not zero,
as it is in (1.2). The choice
of unity as the lower bound of integration dispenses with convergence
problems at that point, while the appearance of the factor $x^{-s}$
instead of the customary $x^{s-1}$ is technically more convenient. Also
it may be compared with the discrete representation
$$
\zeta^{2k}(s) \= \sum_{n=1}^\infty d_{2k}(n)n^{-s}\qquad(\s >
1,\,k\in\NN),
$$
where $d_m(n)$ is the number of ways $n$ may be written as a product of
$m$  factors; $d(n) \equiv d_2(n)$ is the number of divisors of $n$.
Since we have (see [3, Chapter 8])
$$
\int_0^T|\zt|^{2k}\d t \;\ll\; T^{(k+2)/4}\log^{C(k)}T
\qquad(2 \le k \le 6; C(k)\ge0),
$$
it follows that the integral defining ${\Cal Z}_k(s)$ is absolutely
convergent
for $\s > 1$ if $0 \le k \le 2$ and for $\s > (k+2)/4$ if $2 \le k \le 6$.

\medskip
The function $\Z_k(s)$ is a special case of the multiple Dirichlet series
$$
Z(s_1,\cdots,s_{2k},w)
= \int_1^\infty \zeta(s_1+it)\cdots \zeta(s_k+it)\zeta(s_{k+1}-it)
\cdots\zeta(s_{2k}-it)t^{-w}\,{\roman d}t
\leqno(1.3)
$$
considered in a recent work of
A. Diaconu, D. Goldfeld and J. Hoffstein [1].
Analytic properties of this function
may be put to advantage to deal with the important problem of the
analytic continuation of the function $\Z_k(s)$ itself.
It is shown in [1] that (1.3) has meromorphic
continuation (as a function of $2k+1$ complex variables) slightly
beyond the region of absolute convergence, with a polar divisor at
$w=1$. It is also shown that (1.3) satisfies certain quasi-functional
equations, which are used to obtain meromorphic continuation to an
even larger region. Under the assumption that
$$
Z(\hf,\cdots,\hf,w) \;\equiv\; \Z_k(w)
$$
has holomorphic continuation to
the region $\Re{\roman e}\, w \ge 1$ (except for the pole at $w=1$
of order $k^2+1$),
the authors derive the conjecture
on the moments of the zeta-function on the critical line in the form
$$
\int_0^T|\zt|^{2k}\d t \;=\; (c_k + o(1))T\log^{k^2}T\qquad(T\to\infty),\leqno(1.4)
$$
where $k\ge2$ is a fixed integer and
$$
c_k = {a_kg_k\over\G(1+k^2)},\; a_k = \prod_{p}(1-1/p)^{k^2}
\sum_{j=0}^\infty d^2_k(p^j)p^{-j},\; g_k = (k^2)!\prod_{j=0}^{k-1}{j!\over(j+k)!}.
\leqno(1.5)
$$
The formulas (1.4)-(1.5) coincide with the well-known conjecture from Random
Matrix Theory (see e.g., J. Keating and N. Snaith [17]) on the even moments of $|\zt|$.

\medskip
In general one expects, for any fixed $k \in\NN$,
$$
\int_0^T|\zt|^{2k}\d t = TP_{k^2}(\log T) + E_k(T)\leqno(1.6)
$$
to hold (see the author's monograph [4] for an extensive account),
where it is generally assumed that
$$
P_{k^2}(y) = \sum_{j=0}^{k^2}a_{j,k}y^j\leqno(1.7)
$$
is a polynomial in $y$ of degree $k^2$ (the integral in (1.6) is
$\gg_k T\log^{k^2}T$; see e.g., [3, Chapter 9]). The function
$E_k(T)$ is to be considered as the error term in (1.7), namely
one supposes that
$$
E_k(T) \= o(T)\qquad(T \to \infty).\leqno(1.8)
$$
So far (1.6)--(1.8) are known to hold only for $k = 1$ and $k = 2$
(see [3], [4] and [18]).

\medskip
In case (1.6)--(1.8) hold, this may be used to obtain the
analytic continuation of $\Z_k(s)$ to the region $\s\ge 1$ (at least).
Indeed, by using (1.6)--(1.8) we have
$$
\eqalign{&
\Z_k(s) = \int_1^\infty |\zx|^{2k}x^{-s}\d x = \int_1^\infty
x^{-s}\d\left(xP_{k^2}(\log x) + E_k(x)\right)\cr&
= \int_1^\infty (P_{k^2}(\log x) + P'_{k^2}(\log x))x^{-s}\d x
- E_k(1) + s\int_1^\infty E_k(x)x^{-s-1}\d x.\cr}\leqno(1.9)
$$
But for $\R s > 1$ change of variable $\log x = t$ gives
$$
\eqalign{&
\int_1^\infty (P_{k^2}(\log x) + P'_{k^2}(\log x))x^{-s}\d x\cr&
 = \int_1^\infty \left\{\sum_{j=0}^{k^2}a_{j,k}\log^jx +
\sum_{j=0}^{k^2-1}(j+1)a_{j+1,k}\log^jx\right\}x^{-s}\d x\cr&
= \int_0^\infty \left\{\sum_{j=0}^{k^2}a_{j,k}t^j +
\sum_{j=0}^{k^2-1}(j+1)a_{j+1,k}t^j\right\}{\roman e}^{-(s-1)t}\d t\cr&
= {a_{k^2,k}(k^2)!\over(s-1)^{k^2+1}} +
\sum_{j=0}^{k^2-1}(a_{j,k}j! + a_{j+1,k}(j+1)!)(s-1)^{-j-1}.
\cr}\leqno(1.10)
$$
Hence inserting (1.10) in (1.9) and using (1.8)
we obtain  the analytic continuation of $\Z_k(s)$
to the region $\s\ge1$. As we know (see [3], [4], [11] and [19]) that
$$
\int_1^T E_1^2(t)\d t \ll T^{3/2},\qquad
\int_1^T E_2^2(t)\d t \ll T^{2}\log^{22}T,\leqno(1.11)
$$
it follows on applying the Cauchy--Schwarz inequality to the last
integral in (1.9) that (1.8)-(1.10) actually provides
the analytic continuation of $\Z_1(s)$ to the region $\R s > 1/4$, and of
$\Z_2(s)$ to $\R s > 1/2$, but is is actually known that
$\Z_1(s)$ (resp. $\Z_2(s)$) has meromorphic continuation to $\CC$.
For this, see M. Jutila [16] when $k=1$ and Y. Motohashi [19] when $k=2$.

\medskip
The preceding discussion shows one of the several aspects of the
connection between the function $\Z_k(s)$ and power moments of $|\zt|$.
The aim of this paper is to bring forth some results concerning the mean
values of $\Z_k(s)$ and sums of integrals of the form
$$
\sum_{r=1}^R\int_{t_r-G}^{t_r+G}|\zt|^{2k}\d t\quad(T < t_1 < \cdots <t_R< 2T)
\leqno(1.12)
$$
for well-spaced points $t_r$ which satisfy
$t_{r+1} - t_r \ge G\;(r = 1,\ldots, R-1)$,
where $G = G(T)$ is parameter satisfying  $ T^\e \le G \le T^{1-\e}$,
while here and later $\e$ denotes arbitrarily small constants, not
necessarily the same ones at each occurrence.

\medskip
Bounds for sums of the type (1.12) with $k=2$ were
obtained first by H. Iwaniec [15],
who showed that the left-hand side of (1.12) in this case is bounded by
$T^\e(RG + R^{1/2}TG^{-1/2})$ for $T^{1/2} \le G \le T$. Later the author and
Y. Motohashi [12] replaced $T^\e$ by a log-power. In their work [11] the range for
$G$ was relaxed to $\log T \ll G \ll T/\log T$, and the result was generalized.
Further generalizations and results were obtained by the author in [5].

\medskip
One of the applications involving sums of the form (1.12) consists of obtaining
upper bounds for moments of $|\zt|$. Namely one counts (see e.g., Chapter 8 of
[3]) $S$, the number of well-spaced points $\tau_s$ in $[T, \,2T]$ ($\tau_{s+1}
- \tau_s \ge 1)$ such that $|\z(\hf + i\tau_s)| \ge V \,(\ge T^\e)$. Then,
by Theorem 1.2 of [4], it follows that for any fixed $k\in \NN$ we have
$$
V^{2k} \le |\z(\hf + i\tau_s)|^{2k} \ll \log T\int_{\tau_s-1/3}^{\tau_s-1/3}
|\z(\hf + iu)|^{2k}\d u\qquad(s = 1,2,\ldots\,, S), \leqno(1.13)
$$
and one groups integrals on the right hand side of (1.13) into sums of $R$
integrals over intervals $[t_r - G,\,t_r+G]$ with $t_{r+1}-t_r\ge G$ (by
considering separately $r$ with even and odd indices). In this way sums of
the type (1.12) arise, and their estimation leads to estimates for
$\int_0^T|\zt|^{2k}\d t$, which is one of the central problems in the theory
of $\z(s)$.
\bigskip\medskip
\heading {\bf  2.
Statement of results }
  \endheading
We begin with

\bigskip
THEOREM 1. {\it Let $T < t_1 < t _2 < \ldots < t_R < 2T$, $t_{r+1} - t_r \ge G$
for $r = 1, \ldots\,, R-1$. If, for fixed $m,k\in \NN$, we have
$$
\int_T^{2T}\left({1\over G}\int_{t-G}^{t+G}|\zeta(\hf+iu|^{2k}\d u\right)^m\d t
\;\ll_\e\; T^{1+\e}\leqno(2.1)
$$
for $G = G(T) \ge T^{\a_{k,m}}$ and $0\le \a_{k,m}\le 1$, then}
$$
\sum_{r=1}^R\int_{t_r-G}^{t_r+G}|\zt|^{2k}\d t \ll_\e (RG)^{m-1\over m}
T^{{1\over m}+\e}.\leqno(2.2)
$$

\bigskip
The second result, although it could be easily generalized to sums
of the form (1.12), deals with sums of fourth powers. This is because
we have satisfactory results on the mean square of $\Z_k(s)$ so far only for
$k = 1,2.$ The result is

\bigskip
THEOREM 2. {\it Let $T < t_1 < t _2 < \ldots < t_R < 2T$, $t_{r+1} - t_r \ge G$
for $r = 1, \ldots\,, R-1$. Then, for fixed $\hf<\s < 1$, we have}
$$
\sum_{r=1}^R\,\int\limits_{t_r-G}^{t_r+G}|\zt|^{4}\d t
\;\ll_\e\; RG\log^4T + \Bigl(RGT^{2\s-1}
\int\limits_0^{T^{1+\e}G^{-1}}|\Z_2(\s+it)|^2\d t
\Bigr)^{1/2}.\leqno(2.3)
$$

\bigskip
The estimate (2.3) clearly shows the importance of the estimation of $\Z_2(s)$.
Concerning the pointwise estimation of $\Z_2(s)$,
we have (see the author's work [9])
$$
\Z_2(\s + it) \;\ll_\e\; t^{{4\over3}(1-\s)+\e}\qquad(\hf < \s \le 1;\,
t \ge t_0 > 0),\leqno(2.4)
$$
and it was conjectured in [7] that the exponent on
the right-hand side of (2.4) can be replaced by $1/2 - \s$. This
conjecture is very strong, as it was shown in [7] that it implies
$$
\int_0^T|\zt|^8\d t \;\ll_\e\; T^{1+\e},\quad E_2(T) \;\ll_\e\; T^{1/2+\e},
\leqno(2.5)
$$
where $E_2(T)$ (cf. (1.6))
is the error term (see [4], [6], [19]) in the asymptotic formula
for the fourth moment of $|\zt|$. Both estimates in (2.5) are,
up to ``$\e$", known to be best possible.

For the mean square bounds of $\Z_2(s)$ we have the following. It was
proved in by M. Jutila. Y. Motohashi and the author in [14] that
$$
\int_1^T|{\Z}_2(\s + it)|^2\d t \;\ll_\e\;
T^\e\left(T + T^{2-2\s\over1-c}\right) \qquad(\hf < \s \le 1),\leqno(2.6)
$$
and we also have unconditionally
$$
\int_1^T|\Z_2(\s + it)|^2\d t \;\ll\;T^{10-8\s\over3}\log^CT
\qquad(\hf < \s \le 1,\,C > 0).\leqno(2.7)
$$
The constant $c$ appearing in (2.6) is defined by $
E_2(T) \ll_\e T^{c+\e},$ and it is known (see e.g., [4] or [12]) that
$\hf \le c \le {2\over3}$. In (2.6)--(2.7) $\s$ is assumed
to be fixed, as $s = \s+it$ has to stay away from the $\hf$-line
where $\Z_2(s)$ has poles. Lastly, the author [10] proved that,
for ${5\over6} \le \s \le {5\over4}$ we have,
$$
\int_1^{T}|\Z_2(\s+it)|^2\d t \;\ll_\e\; T^{{15-12\s\over5}+\e}.\leqno(2.8)
$$
The lower limit of integration in (2.6)--(2.8) is unity, because of the pole
$s=1$ of $\Z_2(s)$. By taking $c = 2/3$ in (2.6) and using the convexity
of mean values (see e.g., [3, Lemma 8.3]) it follows that
$$
\int_1^T|\Z_2(\s+it)|^2\d t \;\ll_\e\; T^{{7-6\s\over2}+\e}\qquad(\hf < \s
\le \txt{5\over6}).\leqno(2.9)
$$
Note that (2.8) and (2.9) combined  provide the sharpest known bounds
for the mean square of $\Z_2(s)$ in the
whole range $\hf < \s \le {5\over6}$.

\medskip
{\bf Corollary 1.} {\it We have}
$$
\int_0^T|\zt|^{12}\d t\;\ll_\e\; T^{2+\e}.\leqno(2.10)
$$

\medskip
This follows from (2.3) and (2.7) with $\s = 1/2+\e$, giving Iwaniec's bound
$$
\sum_{r=1}^R\,\int\limits_{t_r-G}^{t_r+G}|\zt|^{4}\d t
\ll_\e T^\e(RG + R^{1/2}TG^{-1/2}),
$$
and then taking $k = 2, G = VT^{-\e}$ in conjunction with (1.13). One
immediately obtains $R \ll_\e T^{2+\e}V^{-12}$, and (2.10) follows. This
result (with $\log^{17}T$ replacing $T^\e$) is due to D.R. Heath-Brown [2],
and still represents the strongest bound concerning high moments of $|\zt|$.

\medskip
In obtaining the analytic continuation and bounds for $\Z_2(s)$ in [14],
the authors considered the function
$$
Z_\xi(s) := \int\limits_1 ^\infty J_2(x;x^\xi)x^{-s}\d x,
\; J_k(x;G) := {1\over\sqrt{\pi}G}\int\limits_{-\infty}^\infty
|\z(\hf + ix+iu)|^{2k}{\roman e}^{-(u/G)^2}\d u,\leqno(2.11)
$$
where $k\in \NN,\, 0 < \xi \le 1$, and initially $\R s > 1$. Because of the
smooth Gaussian weight in (2.11) the function $Z_\xi(s)$
is in many aspects less difficult to deal with
than the function $\Z_2(s)$ itself, especially in view of the spectral expansion
of $J_2(x;G)$ obtained by Y. Motohashi (see [18] and [19]). Moreover, by Mellin
inversion and Parseval's formula for Mellin transforms, one can connect bounds
for $Z_\xi(s)$ to the left-hand side of (2.1) when $k=m=2$, and hence
indirectly to power moments of $|\zt|$. Therefore it seems of interest to obtain
bounds for $Z_\xi(s)$, especially if they improve on the existing bounds for
$\Z_2(s)$. In this direction we shall prove in this work a result which is
stronger than the analogous bound (2.4) for $\Z_2(s)$. This is

\bigskip
THEOREM 3. {\it If $\s$ and $\xi$ are fixed, then}
$$
Z_\xi(\s+it) \;\ll_\e\; |t|^{1-\s+\e}
\qquad(\hf < \s \le 1,\,\txt{1\over3} \le \xi \le 1).\leqno(2.12)
$$

\heading
{\bf 3. Proof of Theorem 1 and Theorem 2}
\endheading
We begin with the proof of Theorem 1.
Set $L_k(t,G) = \int_{t-G}^{t+G}|\z(\hf+iu)|^{2k}
\d u.$ Note that if $\mu(\cdot)$ denotes
measure, the bound
$$ \mu\Bigl(t\in[T,\,2T] \,:\,L_k(t,G)\ge
GU\Bigr) \ll_\e T^{1+\e}U^{-m} \qquad(U>0)\leqno(3.1)
$$
follows from the assumption (2.1). We fix $G = G(T)$
and  divide the sum over $r$ in (2.2) into $O(\log T)$ subsums
where $GU < L_k(t_r,G) \le 2UG$. Then, for $U_0\;(\gg1)$ to be
determined later, we have
$$
\eqalign{ &\sum_{r=1}^R L_k(t_r,G) \ll
GRU_0 + \log T\max_{U\ge U_0} \sum_{r,GU < L_k(t_r,G) \le
2GU}L_k(t_r,G)\cr&
\ll GRU_0 + GU\log T\max_{U\ge U_0} \sum_{r,GU <
L_k(t_r,G) \le 2GU}1\cr&
\ll_\e GRU_0 + \log T\max_{U\ge
U_0} T^{1+\e}U^{1-m}\cr& \ll_\e GRU_0 +
T^{1+\e}U_0^{1-m}.\cr}
$$
Here we used the condition that $m\ge1$ and the bound
$$
\sum_{r,L_k(t_r,G)>GU}1 \;\ll_\e\; T^{1+\e}U^{-m}G^{-1}.\leqno(3.2)
$$
To see this, note that if $L_k(t_r,G) > GU$, then
$$
L_k(t,2G) \ge L_k(t_r,G) > GU\qquad({\roman {for}}\;|t-t_r|\le G).
$$
As we can split the sequence of points $\{t_r\}$ into five subsequences,
say $\{t'_r\}$, such that $|t'_{r_1} - t'_{r_2}| \ge 5G$ for
$r_1 \ne r_2$, we see that
$$
G\sum_{r,L_k(t_r,G)>GU}1 \ll
\mu\Bigl(t\in[T,\,2T] \,:\,L_k(t,2G)\ge GU\Bigr),
$$
and (3.2) follows from (3.1). The choice
$$
U_0 \;=\; {\left({T\over RG}\right)}^{1/m}\quad(\,\gg 1\,)
$$
yields
$$ \sum_{r=1}^R
L_k(t_r,G) \;\ll_\e\; T^{1/m+\e}R^{1-1/m}G^{1-1/m},
$$
which is our assertion (2.2).

\medskip
{\bf Corollary 2}. {\it If the hypotheses of Theorem 1 hold, then we have}
$$
\int_0^T|\zt|^{2km}\d t \;\ll_\e\; T^{1+(m-1)\a_{k,m}+\e}.\leqno(3.3)
$$
This follows from (1.13), analogously to Corollary 1.
Observe that ($G = x^\xi, Q_4 = P_4+P'_4$; see (1.6)))
$$
\eqalign{
J_2(x;G) &= {1\over\sqrt{\pi}G}\int_{-\infty}^\infty\Bigl\{
Q_4(\log(x+u) + {\d\over\d u}E_2(x+u)){\roman e}^{-(u/G)^2}\d u\Bigr\}\cr&
= O(\log^4x) + {2\over\sqrt{\pi}G^3}\int_{-\infty}^\infty uE_2(x+u)
{\roman e}^{-(u/G)^2}\d u.\cr}
$$
Hence using the second bound in (1.11) it follows that (2.1), for $k = m = 2$,
holds with $\a_{2,2} = \hf$. By (3.3) this leads to the bound
$$
\int_0^T|\zt|^8\d t \;\ll_\e\;T^{3/2+\e},
$$
which is (up to ``$\e$", see Chapter 8 of [3]) the sharpest one known.

\medskip
We pass now to the proof of Theorem 2. By the inversion formula for the
modified Mellin transform (see Lemma 1 of the author's paper [7]) we have
$$
|\z(\hf+ix)|^4 \= {1\over2\pi i}\int_{(1+\e)}\Z_2(s)x^{s-1}\d s \qquad(x>1).
\leqno(3.4)
$$
In (3.4) we replace the line of integration  by the contour ${\Cal L}$,
consisting of the same straight line from which the segment
$[1+\e-i,\,1+\e+i]$ is removed and replaced by a circular arc
of unit radius, lying
to the left of the line, which passes over the pole $s =1 $ of
the integrand. By the residue theorem we deduce from (3.1) that
$$
|\zx|^4 \= {1\over2\pi i}\int_{\Cal L}{\Cal Z}_2(s)x^{s-1}\d s
+ Q_4(\log x) \qquad(x > 1)
$$
holds, where we have set (cf. (1.6) with $k=2$)
$$
Q_4(\log x) \= P_4(\log x) + P'_4(\log x).
$$
Therefore, for a suitable constant $c$ satisfying $\hf < c < 1$, we have
$$
|\zx|^4 \= {1\over2\pi i}\int_{(c)}{\Cal Z}_2(s)x^{s-1}\d s
+ Q_4(\log x) \qquad(x > 1),\leqno(3.5)
$$
where $\int_{(c)}$ denotes integration over the line $\R s = c$.
Let now $\f_r(x)\;(\ge0)$ be a
smooth function supported in $[t_r-2G,\,t_r+2G]$ such that $\f_r(x) = 1$
when $x\in [t_r-G,\,t_r+G]$, so that
$$
\f_r^{(m)}(x) \;\ll_{r,m}\;G^{-m}\qquad(r = 1,\ldots\,R,\;m = 0,1,2,\ldots\,).
\leqno(3.6)
$$
Analogously as in the proof of Theorem 1, we
can split the sequence $\{t_r\}$ into five subsequences $\{t'_r\}$ such that
that $|t'_{r_1} - t'_{r_2}| \ge 5G$ for $r_1 \ne r_2$. If we multiply (3.5)
by $\f_r(x)$, integrate and sum, we see (writing again $t_r$ for $t'_r$)
that the left hand side of (2.3) is majorized by five sums of the type
$$
\eqalign{&
\sum_{r\le R}\int_{t_r-2G}^{t_r+2G}\f_r(x)|\zx|^4\d x
= O(RG\log^4T) + \cr&
+ \sum_{r\le R}{1\over2\pi i}
\int_{(c)}{\Cal Z}_2(s)\Bigl(\int_{t_r-2G}^{t_r+2G}\f_r(x)x^{s-1}\d x\Bigr)
\,\d s,\cr}\leqno(3.7)
$$
the integrals on the left-hand side of (3.7) being taken over disjoint intervals.
Integrating by parts the integral over $x$ in (3.7) $m$ times, it follows that
it equals
$$
(-1)^m\int_{t_r-2G}^{t_r+2G}\f_r^{(m)}(x){x^{s+m-1}\over
s(s+1)\ldots (s+m-1)}\d x \ll_{r,m} {T^{\s+m-1}\over G^m (1+|t|)^m}.\leqno(3.8)
$$
We can write
$$
\sum_{r\le R}\int_{t_r-2G}^{t_r+2G}\f_r(x)x^{s-1}\d x
= \int_{T/2}^{5T/2}\Phi(x)x^{s-1}\d x,
$$
where $\Phi(x)$ equals $\f_r(x)$ in $[t_r-2G,\,t_r+2G]\,$, and otherwise it
is equal to zero. The bound in (3.8) shows that the portion of the integral in (3.7)
over $s$ for which $|t| \ge T^{1+\e}G^{-1}$ is negligibly
small (i.e., $\ll T^{-A}$ for
any given constant $A>0$), provided that $m = m(\e,A)$ is a sufficiently large
integer.

Thus the left-hand side of (2.3) is, for fixed $\hf < \s < 1$,
$$
\eqalign{&
\ll RG\log^4T + \int_{-T^{1+\e}G^{-1}}^{T^{1+\e}G^{-1}}|\Z_2(\s+it)|
\Biggl|\int_{T/2}^{5T/2}\Phi(x)x^{s-1}\d x\Biggr|\d t\cr&
\ll RG\log^4T + \Biggl(\int_0^{T^{1+\e}G^{-1}}|\Z_2(\s+it)|^2\d t\Biggr)^{1/2}
\Biggl(\int_{T/2}^{5T/2}\Phi^2(x)x^{2\s-1}\d x\Biggr)^{1/2}\cr&
\ll RG\log^4T + \Biggl(\int_0^{T^{1+\e}G^{-1}}|\Z_2(\s+it)|^2\d t\Biggr)^{1/2}
(RGT^{2\s-1})^{1/2},\cr}
$$
which is the assertion of Theorem 2. Here we used, beside the Cauchy-Schwarz
inequality, the estimation
$$
\int_{T/2}^{5T/2}\Phi^2(x)x^{2\s-1}\d x \le
\int_{T/2}^{5T/2}\Phi(x)x^{2\s-1}\d x \ll RGT^{2\s-1},
$$
and the following  (this is Lemma 4 of [7])
\medskip
LEMMA 1. {\it Suppose that $g(x)$ is a real-valued,
integrable function on $[a,b]$, a subinterval
of $[2,\,\infty)$, which is not necessarily finite. Then}
$$
\int\limits_0^{T}\Bigl|\int\limits_a^b g(x)x^{-s}\d x\Bigr|^2\d t
\,\le\, 2\pi\int\limits_a^b g^2(x)x^{1-2\s}\d x \quad(s = \s+it\,,T > 0,\,a<b).
$$
\no This completes the proof of Theorem 2.

\bigskip
\heading \bf
4. The proof of Theorem 3
\endheading
The estimation of $Z_\xi(s)$ was indirectly carried out
in [7] and [9] by the author, in the process of the estimation
of the function $\Z_2(s)$. This function bears resemblance
to the function   $\Z_2(s)$, and it also has a pole of order
five at $s=1$, and infinitely many poles on the line
$\R s = \hf$. For $\Z_2(s)$ Y. Motohashi (see  [19])
showed  that  it has meromorphic
continuation over $\CC$. In the half-plane $\s = \R s >
0$ it has the following singularities: the pole $s = 1$ of order
five, simple poles at $s = {1\over2} \pm i\k_j\,(\k_j =
\sqrt{\lambda_j - {1\over4}})$ and poles at $s = \hf\rho$, where
$\rho$ denotes complex zeros of $\zeta(s)$.  Here as usual
$\,\{\lambda_j = \k_j^2 + {1\over4}\} \,\cup\, \{0\}\,$ is the
discrete spectrum of the non-Euclidean Laplacian acting on
$SL(2,\ZZ)$-automorphic forms (see [19, Chapters 1--3] for a
comprehensive account of  spectral theory and the Hecke
$L$-functions).

\medskip
The estimation of $Z_\xi(s)$ reduces to the estimation of $O(\log t)$
finite integrals of the form
$$
\int_{X/2}^{5X/2}\s(x)J_2(x;x^\xi)x^{-s}\d x,\leqno(4.1)
$$
where ($t> t_0 > 0$ is assumed) as in Section 3 of [7] $t^{1-\e}
\ll X \ll t^A\;(A = A(\s) > 0)$ holds, and $\s(x)\,(\ge0)$ is a smooth
function supported in $[X/2,\,5X/2]$, which equals unity in $[X,2X]$.
For $J_2(x;x^\xi)$ we use Y. Motohashi's spectral decomposition (see [19]),
which we state here as

\medskip
LEMMA 2. {\it If $J_2(x;x^\xi)$ is defined  by} (2.11), {\it then we have}
$$
J_2(T;T^\xi) =I_{2,r}(T,T^{\xi})+I_{2,h}(T,T^{\xi})+I_{2,c}
(T,T^{\xi})+I_{2,d}(T,T^{\xi}).
\leqno(4.2)
$$  {\it
Here $I_{2,r}$ is an explicit main term, the contribution of
$I_{2,h}$ is small,
$$
I_{2,c}(T,T^{\xi })=\pi ^{-1}\int_{-\infty }^{\infty }
{|\zeta (\hf+ir)|^6
\over |\zeta (1+2ir)|^2}\Lambda (r;T,T^{\xi })\d r,
$$
$$
I_{2,d}(T,T^{\xi})=\sum_{j=1}^{\infty}\alpha _j  H_j^3(\hf)\Lambda
(\k_j;T,T^{\xi}),\leqno(4.3)
$$
where
$$
\eqalign{
\Lambda (r;T,T^{\xi}) &= \hf {\roman Re}\Biggl \{\left (1+
{i\over \sinh \pi r} \right )\Xi (ir;T,T^{\xi}) \cr&+ \left (1-
{i\over \sinh \pi r} \right )\Xi (-ir;T,T^{\xi}) \Biggr \}
\quad(r\in\RR)\cr }
\leqno(4.4)
$$
with
$$
\eqalign{
\Xi (ir;T,T^{\xi }) &={\G^2 ({1\over2}+ir)\over{\G (1+2ir)}}\int_0^{\infty }
(1+y)^{-{1\over2}+iT}y^{-{1\over2}+ir}\cr
& \times \exp\left(-\txt{1\over4} T^{2\xi}\log ^2(1+y)\right )
F(\hf +ir, \hf + ir;1+2ir;-y) \d y,\cr}
\leqno(4.5)
$$
and $F$ is the hypergeometric function}.

\medskip
The contribution of the main term $I_{2,r}$ in (4.2)
(of order $\ll \log^4T$) to (4.1) is small if one uses integration by parts
and $\s^{(m)}(x) \ll_m X^{-m}\;(m\ge0)$.
The same is true of the contribution of the
continuous spectrum $I_{2,c}$, if one uses the bounds for
$\Xi (ir;T,T^{\xi })$ in Chapter 5 of [19]. The main contribution comes
from $I_{2,d}$ (the ``discrete spectrum") in (4.2), and the problem reduces
to the asymptotic evaluation of the functions $\Lambda (r;T,T^{\xi})$
and $\Xi (ir;T,T^{\xi })$ in (4.4) and (4.5), respectively. This task
was carried out in detail in the recent work of A. Ivi\'c--Y. Motohashi
[13]. In particular, we invoke the discussion in Section 5 of this paper.
The major contribution to (4.1) of $\Xi (ir;T,T^{\xi })$, by equation
(5.14) of [13] turns out to be a multiple of
$$
\int_{X/2}^{5X/2}\s(x)x^{-s}\Bigl(\sum_{\k_j\le X^{1-\xi}\log X}\a_j\H\,
 I_\xi(x,\k_j)\Bigr)\d x,
\leqno(4.6)
$$
where, for any fixed integer $N$, $G = x^\xi$ and effectively computable
constants $c_j$,
$$
\eqalign{&I_\xi(x,\k_j) :=
x^{-1/2}\k_j^{-1/2}
\exp\left\{-{\txt{1\over4}}G^2\log^2(1+y_0)
+i{\Cal F}(y_0)-i\k_j\log4\right\}\cr&
= x^{-1/2}\k_j^{-1/2}\times\cr&
\exp\left\{-{\txt{1\over4}}G^2\log^2(1+y_0)
+ i\k_j\log\left({\k_j\over4{\roman e}x}\right) + i\sum_{j=3}^Nc_j\k_j^jx^{1-j}
   + O_N(\k_j^{N+1}x^{-N})\right\}.\cr}\leqno(4.7)
$$
This is understood in the following sense: the remaining terms in
the evaluation of the relevant expression are either negligible (smaller
than $X^{-A}$ for any constant $A>0$), or similar in nature to
(4.7) (meaning that the oscillating exponential factor is the same,
which is crucial), only of the lower order of magnitude than
the corresponding terms in (4.7).
The remaining notation is as follows. We have
$$
y_0 = {\k_j\over x}\left(\sqrt{1 + {\k_j^2\over4x^2}} +
{\k_j\over2x}\right), \leqno(4.8)
$$
so that $y_0 \sim \k_j/x$ as $x\to\infty$ in the relevant range. Moreover, we have
$$
{\Cal F}(y_0) = \k_j\log y_0 - 2\k_j\log\left({1+\sqrt{1+y_0}\over2}\right)
- x\log(1+y_0).\leqno(4.9)
$$
The term $\exp\Bigl(O_N(\k_j^{N+1}x^{-N})\Bigr)$ in (4.7) is expanded into
a power series. If we take $N$ sufficiently large, then only the first term
unity will make a non-negligible contribution. Hence instead of (4.6) we
need to estimate
$$
\sum_{\k_j\le X^{1-\xi}\log X}\a_j\k_j^{-1/2}\H\int_{X/2}^{5X/2}x^{-1/2-i\k_j-s}
\Bigl(\s(x) L_\xi(x,\k_j)\Bigr)\d x,
\leqno(4.10)
$$
say, where
$$
L_\xi(x,\k_j):= \exp\left\{-{\txt{1\over4}}G^2\log^2(1+y_0)
+ i\k_j\log\left({\k_j\over4{\roman e}}\right) +
i\sum_{\ell=3}^Nc_\ell\k_j^\ell x^{1-\ell}\right\}.\leqno(4.11)
$$
We integrate (4.10) many times by parts, using (4.11) and the facts that
$\s(X/2) = \s(5X/2) = 0$ and $\s^{(m)}(x) \ll_m X^{-m}$ for $m \ge 0$.
Thus, since the integral of $x^{-1/2-i\k_j-s}$ is
$$
{x^{1/2-i\k_j-s}\over  {1\over2}-i\k_j-s}\qquad(\s >\hf, \;\hf-i\k_j-s \ne 0)
$$
and
$$
(\s(x) I_1)' \ll X^{-1} + \k_j^3X^{-3} \ll X^{-1} + X^{-3\xi}\log^3X
\ll X^{-1}\log^3X
$$
for $\k_j \le X^{1-\xi}\log X$ and $\xi \ge 1/3$ (which is our
assumption for this reason), it follows that
only the values of $\k_j$ for which $|\k_j - t| \ll_\e t^\e$ will make
a non-negligible contribution. To complete the proof we need now
(see the author's paper [8]) the bound contained in

\medskip
LEMMA 3. {\it We have}
$$
\sum_{K-G\le\k_j\le K+G} \a_j\H \;\ll_\e\; GK^{1+\e}
\quad(K^{\e}  \;\le\; G \;\le \; K).\leqno(4.12)
$$

\medskip
Thus we are left with the contribution which is, by (4.12),
$$
\ll_\e\; t^{-1/2}\sum_{|t-\k_j|<t^\e}\a_j\H X^{1/2-\s}
\;\ll_\e\; t^{1/2+\e}X^{1/2-\s}.
$$
Since $\s > \hf$, the last expression is $\ll_\e t^{1-\s+\e}$ in the
relevant range $X \gg t^{1-\e}$, and (2.12) follows. Our result is
certainly not optimal, since by using (4.12) we have ignored the
exponential factor in $L_\xi(x,\k_j)$ in (4.11) and the factor $x^{-i\k_j}$ in
(4.10). On the other hand, there do not exist yet non-trivial estimates for
exponential sums with $\a_j\H$, which vitiates our efforts to improve
on (2.12).

\vfill
\eject

\Refs


\item{[1]} A. Diaconu, D. Goldfeld and J. Hoffstein,
Multiple Dirichlet series and moments of zeta and $L$-functions,
{\it Compositio Math.} {\bf139}, No. 3. 297-360(2003).

\item{[2]} D.R. Heath-Brown, The twelfth power moment of the
Riemann zeta-function, {\it Quart. J. Math. (Oxford)} {\bf29}(1978), 443-462.

\item {[3]} A. Ivi\'c,  The Riemann zeta-function, {\it John Wiley
and Sons}, New York, 1985 (2nd ed. {\it Dover}, Mineola N.Y., 2003).

\item {[4]} A. Ivi\'c,  Mean values of the Riemann zeta-function,
LN's {\bf 82}, {\it Tata Institute of Fundamental Research},
Bombay, 1991 (distr. by Springer Verlag, Berlin etc.).

\item{[5]} A. Ivi\'c, Power moments of the Riemann
zeta-function over short intervals, {\it Archiv
     Math.} {\bf62} (1994), 418-424.

\item{[6]} A. Ivi\'c, On the error term for the fourth moment of the
Riemann zeta-function, {\it J. London Math. Soc.}
{\bf60}(2)(1999), 21-32.

\item{ [7]}  A. Ivi\'c, On some conjectures and results
for the Riemann zeta-function
and Hecke series, {\it Acta Arith.}  {\bf109}(2001), 115-145.

\item{[8]} A. Ivi\'c, On sums of Hecke series in short intervals,
{\it J. de Th\'eorie des Nombres Bordeaux} {\bf13}(2001), 1-16.

\item{ [9]}  A. Ivi\'c, On the estimation of ${\Cal Z}_2(s)$,
in ``{\it Anal. Probab. Methods
Number Theory}" (eds. A. Dubickas et al.), TEV, Vilnius, 2002,  83-98.

\item{[10]} A. Ivi\'c, On the estimation of some Mellin transforms
connected with the fourth moment of $|\zt|$, to appear in the {\it Proc.
Conf. ``Elementare und Analytische Zahlentheorie"}, Mainz, May 2004
(ed. W. Schwarz), {\tt ArXiv:math.NT/0404524}.

\item{[11]} A. Ivi\'c and Y. Motohashi, The mean square of the error term for the
    fourth moment of the zeta-function, {\it Proc. London Math. Soc.} (3)
    {\bf69} (1994), 309-329.

\item{[12]} A. Ivi\'c and Y. Motohashi, On the fourth power moment of the
    Riemann zeta-function, {\it J. Number Theory}  {\bf51} (1995), 16-45.

\item{[13]} A. Ivi\'c and Y. Motohashi, The Moments of the Riemann Zeta-Function.
 Part I: The fourth moment off the critical line, Functiones et Approximation
 (subm.), {\tt ArXiv.math.NT/0408022}.

\item{[14]} A. Ivi\'c, M. Jutila and Y. Motohashi, The Mellin transform of
power moments of the zeta-function, {\it Acta Arith.}
 {\bf95}(2000), 305-342.

\item{[15]} H. Iwaniec, Fourier coefficients of cusp forms and the Riemann
zeta-function, Expos\'e No. {\bf18}, {\it S\'eminaire de Th\'eorie des Nombres,
Universit\'e Bordeaux. }1979/80.

\item {[16]}  M. Jutila, The Mellin transform of the square of Riemann's
zeta-function, {\it Periodica Math. Hung.} {\bf42}(2001), 179-190.

\item{[17]} J.P. Keating and  N.C. Snaith, Random matrix theory and
$\zeta(\hf+it)$, {\it Comm. Math. Phys.} {\bf214}(2000), 57-89.

\item{ [18]} Y. Motohashi,   An explicit formula for the fourth power
mean of the Riemann zeta-function, {\it Acta Math. }{\bf 170}(1993), 181-220.

\item {[19]} Y. Motohashi,  Spectral theory of the Riemann
zeta-function, {\it Cambridge University Press}, Cambridge, 1997.

\bigskip

Aleksandar Ivi\'c

Katedra Matematike RGF-a

Universitet u Beogradu

\DJ u\v sina 7, 11000 Beograd

Serbia and Montenegro, \tt ivic\@rgf.bg.ac.yu

\endRefs
\bye